\documentclass{article}
\usepackage{amsmath,amssymb}

\newtheorem{theorem}{Theorem}
\newtheorem{lemma}{Lemma}

\newtheorem{example}{Example}

\begin{document}

\title{A Bijection for Partitions with Initial Repetitions}

\author{William J. Keith \\
University of Lisbon, Av. Prof. Gama Pinto 2, Lisbon, Portugal\\
wjk150@cii.fc.ul.pt}

\maketitle

\begin{abstract} A theorem of Andrews equates partitions in which no part is repeated more than $2k-1$ times to partitions in which, if $j$ appears at least $k$ times, all parts less than $j$ also do so.  This paper proves the theorem bijectively, with some of the generalizations that usually arise from such proofs.
\end{abstract}

\section{Introduction}

In \cite{Andrews1}, George Andrews proved the $q$-series identity

\begin{multline}\label{AndrewsIdent}\sum_{n=0}^{\infty} \frac{q^{k\cdot 1+k\cdot 2+\cdots +k\cdot n}}{(1-q)(1-q^2)\cdots(1-q^n)} \prod_{j=n+1}^{\infty} (1+q^j+q^{2j}+\cdots+q^{(k-1)j}) \\ = \prod_{j=1}^{\infty} (1+q^j+q^{2j} + \cdots + q^{(2k-1)j})
\end{multline}

\noindent which can be interpreted combinatorially as stating that partitions of $n$ in which no part is repeated more than $2k-1$ times (the right-hand side) are equinumerous with partitions of $n$ in which, if any part $j$ appears at least $k$ times as a part, then any part of size less than $j$ also appears at least $k$ times as a part.  These latter he named \emph{partitions with initial $k$-repetitions}.  He requested a proof via bijection of the two sets and kindly noted that before publication the current author (his student at the time) was able to provide.  That proof has not previously appeared in print and the details were recently requested, so it is offered here for the literature.

As is typical of bijective proofs, the map can modestly generalize the theorem by careful observation of the characteristics of the starting and ending partitions.

\section{The Bijection}

We assume familiarity with the common vocabulary of partitions; see \cite{Andrews2} for the standard reference.

A somewhat less-common construction for a given partition is its \emph{$k$-modular diagram}.  In this construction, instead of displaying a part $\lambda_i$ as a column of $\lambda_i$ units, we display $\lambda_i$ by writing its residue modulo $k$, if nonzero, at the top of a column, followed by $\vert \lambda_i \vert_k$ instances of $k$.  Here $\vert \lambda_i \vert_k$ is the largest number of multiples of $k$, i.e. $\lambda_i = \vert \lambda_i \vert_k \cdot k + r$, with $r$ the least nonnegative residue of $\lambda_i$ mod $k$.

\begin{example}\label{kModDiag} Let $\lambda = (29,27,25,21,17,8,8,5,4,1)$.  The 5-modular diagram of $\lambda$ is \[ \, \begin{matrix} 4 & 2 & 5 & 1 & 2 & 3 & 3 & 5 & 4 & 1 \\ 5 & 5 & 5 & 5 & 5 & 5 & 5 & & & \\ 5 & 5 & 5 & 5 & 5 & & & & & \\  5 & 5 & 5 & 5 & 5 & & & & & \\ 5 & 5 & 5 & 5 & & & & & & \\ 5 & 5 & & & & & & & &  \end{matrix} \text{.} \]
\end{example}

The Ferrers diagram of a partition may essentially be regarded as its 1-modular diagram.  The 2-modular diagrams, due to the fact that all nonzero residues are 1, often have convenient combinatorial properties for proofs of partition facts involving partitions into odd or distinct parts.  The notion of a $k$-strip, defined below and central to the current proof, is found in Stockhofe \cite{Stockhofe}, which makes extensive use of $k$-modular diagrams to construct the group of all possible bijections on the set of partitions.  That thesis is in German; an unofficial translation can be found as an appendix to the current author's doctoral thesis, \cite{KThesis}, in which the latter half of the main material employs other concepts therefrom.

If $\lambda_i - \lambda_{i+1} \geq k$, then $\vert \lambda_i \vert_k - \vert \lambda_{i+1} \vert_k \geq 1$.  In that case, it is possible to subtract $k$ from all parts $\lambda_1$ through $\lambda_i$ and still be left with a partition, now of $n-ki$.  Visually, this corresponds to removing a row of $k$-units from the $k$-modular diagram of $\lambda$.  It is possible to have $\vert \lambda_i \vert_k - \vert \lambda_{i+1} \vert_k \geq 1$ with $\lambda_i - \lambda_{i+1} < k$, such as occurs between parts $\lambda_4$ and $\lambda_5$ in Example \ref{kModDiag}, but in this case it is not possible to remove such a row: $\lambda_4 - 5 = 16 < \lambda_5 = 17$.

When it is possible to remove such a row of $k$-units, we will refer to the removable units as a \emph{$k$-strip} of length $i$.  If $\lambda_i - \lambda_{i+1} \geq 2k$, it is possible to have multiple removable $k$-strips of the same length, and if $\lambda$ has $m$ parts, it is possible to have removable $k$-strips of length $m$, treating $\lambda_{m+1} = 0$ by convention.

In Example \ref{kModDiag}, there is one removable $k$-strip, which is of length 5.

If we remove all $k$-strips from a partition (the order in which strips are removed does not affect the number or sizes of strips obtained), we can break down any partition uniquely into an ordered pair of partitions $\lambda \rightarrow (\pi, \delta)$.  Here $\delta$ is a partition into multiples of $k$ which records the amounts removed in each strip, and $\pi$ is what Stockhofe terms \emph{$k$-flat}: its first part is less than $k$ and all differences between adjacent parts are less than $k$.

In Example \ref{kModDiag}, we break down $\lambda$ as $((24,22,20,16,12,8,8,5,4,1), (25))$.

\subsection{Partitions with initial $k$-repetitions}

Recall that conjugation visualized as reflection applies to the original Ferrers diagram of a partition, not its $k$-modular diagram; nonzero residues can only appear in the top row.  However, if there are no nonzero residues, as in $\delta$, reflection about the main diagonal, while not standard conjugation, does reliably result in a unique partition.  In fact, we have essentially done this in constructing $\delta$, a list of sizes of parts, from the extraction of $k$-units by rows out of $\lambda$.

The conjugates of $k$-flat partitions are partitions in which no part appears $k$ or more times.  The set of partitions of $n$ in which no part appears more than $2k-1$ times is mapped by conjugation to the set of partitions of $n$ in which no difference is $2k$ or more, and the smallest part is $2k-1$ or less.  This is exactly the set of partitions in which any $k$-strips are of distinct lengths.

Let $\lambda$ be a partition in which no part appears more than $2k-1$ times.  The bijection which proves Andrews' theorem is as follows:

\begin{itemize}
\item Conjugate $\lambda$, obtaining $\lambda^{\prime}$.
\item Remove all $k$-strips from $\lambda^{\prime}$, obtaining $(\pi, \delta)$.  The partition $\pi$ is $k$-flat, and $\delta$ has distinct parts.
\item Construct $\alpha = \pi + \delta$ by vector addition, i.e. $\alpha_1 = \pi_1 + \delta_1$, $\alpha_2 = \pi_2 + \delta_2$, etc.  Assume the convention that the shorter of $\pi$ or $\delta$ is filled out with parts of size 0.
\item Conjugate $\alpha$.
\end{itemize}

\begin{example} If we begin with the partition \[ (10,9,9,9,8,7,7,7,5,5,5,5,5,5,5,5,5,4,4,4,4,3,3,3,3,2,2,1,1) \]

\noindent then $\lambda^{\prime}$ is the partition of Example \ref{kModDiag}, and so $\pi = (24,22,20,16,12,8,8,5,4,1)$, $\delta = (25)$ as before.  Addition gives us $\alpha = (49,22,20,16,12,8,8,5,4,1)$, so $\alpha^{\prime}$ is our original partition with five fewer parts of size 5 and twenty-five more parts of size 1, i.e. \[ \alpha^{\prime} = (10,9,9,9,8,7,7,7,5,5,5,5,4,4,4,4,3,3,3,3,2,2,1,<\text{25 more} > ,1) \, \text{.} \]
\end{example}

\begin{theorem} The image under this map of a partition with parts repeated at most $2k-1$ times is a unique partition with initial $k$-repetitions. \end{theorem}

\noindent \textbf{Proof:} Since $\pi$ is $k$-flat, $\pi_1 - \pi_2 < k$.  This is precisely the number of times that a part of size 1 appears in $\pi^{\prime}$.  Before any other change is made to $\pi$, no part of any size appears $k$ or more times in $\alpha^{\prime}$.

Suppose $\delta$ has at least 1 part.  When we add $\delta_1$ to $\pi_1$, the number of 1s in $\alpha^{\prime}$ increases by a multiple of k.  However, no larger part in $\alpha^{\prime}$ yet appears $k$ or more times.  So far, the conjugate of $\alpha$ is a partition with initial $k$-repetitions; specifically, only 1 appears more than $k$ times.

Suppose $\delta$ has at least 2 parts.  When we add $\delta_2$ to $\pi_2$, the number of 2s in $\alpha^{\prime}$ increases by $\delta_2$, a multiple of $k$, and the number of 1s in $\alpha^{\prime}$ decreases by $\delta_2$.  However, since $\pi_1 \geq \pi_2$, $\delta_1 > \delta_2$ and $k \vert \delta_i$, the number of parts of size 1, $\pi_1 - \pi_2 + \delta_1 - \delta_2$, is still at least $k$.  Parts of size 2 and 1 in $\alpha^{\prime}$ are both repeated at least $k$ times, but no larger part appears $k$ or more times.

If $\delta$ has at least 3 parts, the addition of $\delta_3$ similarly affects the number of repetitions of parts of size 3 and 2, but neither this nor any further part affects the number of parts of size 1.  This continues up to the number of parts of $\delta$.  (If the number of parts of $\delta$ exceeds the number of parts of $\pi$, then $\alpha^{\prime}$ has largest parts repeated some multiple of $k$ times.)  Thus, $\alpha^{\prime}$ is a partition with initial $k$-repetitions. 

The map is one-to-one since the $k$-flat portion of a partition is uniquely defined by its sequence of residues mod k, and with $\pi$ uniquely defined, $\delta$ is unique.  \hfill $\Box$

\vspace{0.1in}

For the map to re-prove the theorem, it must be a bijection:

\begin{lemma} The inverse of this map is well-defined.\end{lemma}

\noindent \textbf{Proof:} The ideas are all the same.  From a partition in which $j$ and all lower parts appear at least $k$ times but no higher part does so, remove $j$ rows of distinct multiples of $k$, reducing the number of repetitions below $k$ for all of these small parts.  Conjugate the remaining partition to obtain a $k$-flat partition, and add the removed rows back in horizontally, as $k$-strips of distinct lengths.  The resulting partition will have differences less than $2k$. \hfill $\Box$

\vspace{0.1in}

Depending on an investigator's interest, various characteristics can be isolated from this map.  For example, we can finitize the sum on the left-hand side of identity (\ref{AndrewsIdent}) by observing only those partitions with initial $k$-repetitions in which parts no larger than $m$ are repeated $k$ or more times.  These will map to those partitions with differences no larger than $2k-1$ from which no more than $m$ $k$-strips can be extracted:

\begin{multline} \sum_{n=0}^{m} \frac{q^{k\cdot 1+k\cdot 2+\cdots +k\cdot n}}{(1-q)(1-q^2)\cdots(1-q^n)} \prod_{j=n+1}^{\infty} (1+q^j+q^{2j}+\cdots+q^{(k-1)j}) \\ = \prod_{j=1}^{\infty} \frac{1-q^{jk}}{1-q^j} \cdot \sum_{n=0}^{m} \frac{q^{kn(n+1)/2}}{(1-q^k)(1-q^{2k})\cdots(1-q^{nk})}
\end{multline}

Suppose we require that, if any $j$ appear at least $k$ times, $j-1$ and all lower parts appear at least $2k$ times.  Then the $k$-strips removed will differ in length by at least 2, a condition equivalent to the first Rogers-Ramanujan theorem:

\begin{multline} \sum_{n=0}^{\infty} \frac{q^{k\cdot n+2k\cdot (n-1)+\cdots +2k \cdot 1}}{(1-q)(1-q^2)\cdots(1-q^n)} \prod_{j=n+1}^{\infty} (1+q^j+q^{2j}+\cdots+q^{(k-1)j}) \\ = \prod_{j=1}^{\infty} \frac{1-q^{jk}}{1-q^j} \cdot \sum_{n=0}^{\infty} \frac{q^{kn^2}}{(1-q^k)(1-q^{2k})\cdots(1-q^{nk})} \\ = \prod_{j=1}^{\infty} \frac{1-q^{jk}}{1-q^j} \cdot \prod_{n=0}^{\infty} \frac{1}{(1-q^{k(5n-4)})(1-q^{k(5n-1)})} \\ =\prod_{j=1}^{\infty} \frac{(1-q^{k(5j)})(1-q^{k(5j+2)})(1-q^{k(5j+3)})}{1-q^j}
\end{multline}

By setting various conditions on the starting or ending partitions, a whole library of similar manipulations can be invoked.

\end{document}